\documentclass[11pt]{article}

\usepackage{amsmath,amsthm,amsfonts,amssymb,amscd}

 \usepackage[portuguese]{babel}
 \usepackage[latin1]{inputenc}

\newtheorem{teorema}{Teorema}

\newcommand{\R}{\mathbf R}

\newcommand{\fdd }{ \hfill $\Box$ }

\begin{document}

\title{
Equações diferenciais pantográficas: retardo proporcional
}

\date{\today\footnote{{\jobname}.tex}}


\author{Paulo Ruffino \\  Departamento de Matemática \\ Universidade 
Estadual de Campinas \\ ruffino@ime.unicamp.br}

\maketitle

\section{Introdução e Definição}

A intenção dessas notas é mostrar aos leitores um tipo curioso de 
equações
diferenciais ordinárias com retardo, isto é, equações diferenciais cujas derivadas 
dependem, de uma maneira peculiar, do passado da solução. São as chamadas 
equações com retardo proporcional, ou equações pantográficas.

Como motivação geral para equações com retardo, não faltam 
exemplos de fenômenos na natureza que são modelados desta maneira. Se uma 
estrela tem sua massa 
variando com o tempo, dadas as distâncias astronômicas envolvidas, os planetas que orbitam ao redor dela só sentem essa 
variação algum tempo depois. Pensando no nosso Sol, por exemplo: uma 
hipotética explosão demoraria oito minutos para ser sentida no nosso planeta. 
Modelos  
biológicos, econômicos, bioquímicos, hidráulicos e vários outros sistemas, 
intrinsecamente tem um 
comportamento 
com retardo. Até mesmo circuitos eletrônicos apresentam atrasos de 
nanosegundos que em altas frequëncias precisam ser considerados. Esse tipo 
de raciocínio pode ser levado ao extremo a ponto de dizermos que 
todo 
fenômeno da natureza, mesmo nos problemas clássicos, apresenta algum retardo: uns 
mais perceptíveis do que outros.

Tipicamente, nos modelos de equações com retardo existe uma função (funcional,
linear ou não) que associa para cada trajetória em um intervalo de tempo no ``passado
imediato'' dessa solução, uma velocidade da solução no ponto em
consideração. 

%

O tipo de equação diferencial com retardo mais conhecida é a com retardo constante, isto é, em que a derivada do sistema
no instante $t$ depende somente de onde o sistema estava no instante $t-r$, 
para um certo tempo de retardo fixo $r>0$. Mais precisamente: seja $\eta: 
[-r, 0) \rightarrow \R$ uma função contínua em $\R$ (a
condição inicial) e considere
o seguinte problema: encontrar uma função contínua $x: [-r,
\infty) \rightarrow \R$ que satisfaça 

\begin{equation} \label{Eq: retardo constante}
x'(t)= F(x(t-r))
\end{equation}
para $t>0$ com $x(t)=\eta (t)$ para $t\in [-r, 0]$, onde $F$ é uma função 
(campo de vetores)
contínua. Equações
deste tipo já são relativamente bem conhecidas, figuram em textos clássicos da 
área, ver por exemplo J. Hale \cite{Hale} e as referências ali contidas. 
Variações dessas equações incluem 
sua extensão para variedades diferenciáveis (via transporte paralelo ou não, 
por exemplo Oliva \cite{Oliva}) e adição de perturbações estocásticas, ver por 
exemplo \cite{Langevin-Oliva},  \cite{Salah},  entre outros trabalhos, desses e de outros autores.


Neste texto queremos chamar a atenção dos leitores para um outro tipo
especial de equações 
com
retardo, relativamente pouco conhecido, onde o parâmetro no tempo que
define o campo de vetor que dirige sua solução tem um retardo crescente, que é
proporcional ao tempo. 
Fixamos a proporção do retardo $q \in (0,1)$ e um ponto base $a\geq 0$. 
Considere uma função inicial $\eta: [qa, 
a]\rightarrow \R$ contínua. Dada uma 
função (campo de vetores) $F: \R \rightarrow \R$, 
considere a equação diferencial com retardo:
\begin{equation} \label{Eq: retardo pantografico}
x'(t)= F(x(qt)). 
\end{equation}
A solução é uma função contínua $x(t)$ definida em algum 
sub-intervalo de  
$(qa, +\infty)$ de tal maneira que para $t\in [qa, a]$, temos $x(t)= 
\eta(t)$. Equações deste tipo são chamadas por alguns autores 
de {\it equações diferenciais pantográficas}, ver por exemplo \cite{Iserles-Liu-1994}, \cite{Iserles-Liu} e \cite{Kato-McLeod}.

A equação pantográfica modela por exemplo a situação de duas partículas se 
movendo na reta, uma mais lenta  que a outra por um fator $q \in (0,1)$, onde a 
partícula que está a frente, mais rápida, obedece comandos de velocidade dados 
em função da posição da segunda partícula, mais lenta. Outras situações onde 
as equações pantográficas podem aparecer é em sistema com duas escalas de 
tempo, por exemplo nas estruturas onde se estudam os chamados ``princípios de 
médias''. Na próxima seção 
mostraremos propriedades básicas interessantes dessas equações. 



\section{Propriedades}

\subsection{Existência e unicidade de solução}

A primeira propriedade que mostraremos é que, diferente das EDO's clássicas, 
basta que a função $F$ seja contínua para existir uma única solução. 
Além disso, com essa hipótese, a solução não explode em tempo finito, i.e. 
está definida para todo $t\geq qa$.

\begin{teorema} Considere a equação com 
retardo pantográfico (\ref{Eq: retardo pantografico}), 
onde a função  $F$ é contínua.
Dada uma condição inicial contínua (integrável) $\eta: [qa, a] \rightarrow \R$, 
com $a\geq 0$ e 
$q\in (0,1)$
existe uma única 
solução $x(t)$ para essa equação no intervalo $[qa, +\infty)$ tal que $x(t)= 
\eta (t)$ se $t \in [qa, a]$. 
\end{teorema}

\noindent {\it Demonstração:} A solução pode facilmente ser 
construída. Mostraremos por indução que para todo $k \in \{-1, 0, 1, 
\ldots \}$ a solução existe nos intervalos $[aq^{-k}, aq^{-(k+1)}]$. De fato, 
 no intervalo correspondente a $k=-1$ temos a função inicial $\eta(t)$. 
Para $k\geq 0$, nossa hipótese de indução é que $x(t)$ esteja 
definido em $[aq^{-(k-1)},  aq^{-k}]$, então para  $t\in [aq^{-k}, 
aq^{-(k+1)}]$, 
defina
\[
 x(t)= x (aq^{-k}) + \int_{aq^{-k}}^t \ F(x(qs))\ ds.
\]
O integrando está bem definido neste intervalo. Além disso, a continuidade da 
$F$ garante também que $x(t)$ está bem definido (e não vai para infinito) em 
$t\in [aq^{-k}, aq^{-(k+1)}]$.

Para verificarmos a unicidade, suponha que $x(t)$ e $\tilde{x}(t)$ sejam 
soluções para a mesma função inicial $\eta(t)$ no intervalo $[qa, a]$. Se o 
conjunto $D=\{t\geq qa: x(t) \neq \tilde{x}(t) \}$ for diferente de 
vazio, então $D$ possui um ínfimo $t_0 \in \R$. Por continuidade de $x(t)$ e 
$\tilde{x}(t)$, existe um $t\in (t_0,  t_o q^{-1})\cap D$. Então, 
para esse valor de $t$:
\[
 x(t) = x(t_0) + \int_{t_0}^t F(x(qs)) \ ds = \tilde{x}(t_0) + \int_{t_0}^t 
F(\tilde{x}(qs)) \ ds = \tilde{x}(t)
\]
já que $\tilde{x}(qs) = {x}(qs)$ para $s\in [t_0, t_0 q^{-1}]$, que contém o  
intervalo de integração. Essa contradição mostra que o conjunto $D$ é vazio.

\fdd

\subsection{Condição inicial degenerada para um ponto}

Diferente das equações com retardo constante (\ref{Eq: 
retardo 
constante}), curiosamente nossas equações com retardo pantográfico 
(\ref{Eq: 
retardo pantografico}) podem ter suas 
condições (funções) iniciais degeneradas num único ponto: em $t=0$ 
 defina $x(0)= x_0\in \R$, como nas EDO's clássicas. Isso 
corresponde 
a tomar o parâmetro $a=0$ no teorema acima e $\eta(a)=x(0)$.  

A título de exemplo, vamos descrever a solução de uma equação diferencial 
pantográfica linear em termos de sua série de Taylor centrada em $a=0$. 
Considere então a equação com condição inicial:
\begin{equation} \label{Eq: retardo linear}
\left\{ \begin{array}{ccl}
          x' & = &\lambda x(qt), \\
            x(0)&= & x_0 .
\end{array} \right. 
\end{equation}
Então uma solução fundamental $\varphi_t: \R \rightarrow \R$ tal que 
$\varphi'(t)= \lambda \varphi (q t)$ é infinitamente derivável e sua 
$n$-ésima derivada satisfaz 

\[
 \varphi^{(n)}(t)= \lambda^n q^{1+2+ \ldots + (n-1)} \varphi(qt)
\]
Assumindo $\varphi(0)=1$ temos então 
\[
 \varphi^{(n)}(0)= \lambda^n q^{\frac{n(n-1)}{2}}.
\]
Portanto a fórmula de Taylor nos dá a solução fundamental na forma
\begin{equation} \label{solucao 
fundamental}
 \varphi (t)= 1 + (\lambda t) + \frac{q(\lambda t)^2}{2!} + \frac{q^3 (\lambda 
t)^3}{3!} + \ldots + 
\frac{q^{\frac{n(n-1)}{2}} (\lambda t)^n}{n!} + \ldots 
\end{equation}
que converge absolutamente para todo $t \in \R$. Para $\lambda >0$ temos $x(t)$ 
crescente com $\lim_{t\rightarrow \infty} x(t)= \infty$. Para $\lambda<0$ 
temos $x(t)$ decrescente em uma vizinhança do zero e cruza o zero (ver Teorema 
3 e item 1 dos Comentários Finais).

Para ver uma fórmula de Taylor numa situação mais geral, incluindo 
um termo extra linear sem retardo veja por exemplo \cite[Eq. 2.2]{Iserles-Liu}.

\subsection{Propriedade de reconstrução para $t<qa$ (voltando no tempo)}

Tipicamente, nas equações com retardo, não faz sentido perguntarmos por uma 
solução com o tempo antes da função inicial, isto é, para $t<qa$, já que 
não temos
informação sobre a derivada nesse intervalo. No 
entanto, em situações especiais podemos de fato estender essa solução para 
algum intervalo antes de $qa$. Essa é a chamada propriedade de reconstrução das 
equações com retardo. Para isso, suponha que 
\begin{description}
 \item[1)] a função
inicial $\eta: [qa, a] \rightarrow \R $ seja de classe $C^1$ em $[qa, a]$,
 \item[2)] a derivada 
lateral de $\eta$ em $a$ satisfaz $\eta'(a)= F(\eta (qa))$,
 \item[3)] a função de coeficientes (campo de vetores) $F: \R \rightarrow \R$ 
seja inversível.
\end{description}
Verifique que temos então uma solução contínua $x(t)$ em $[q^2 a, \infty)$ tal 
que em 
$t\in [q^2a, qa]$ temos
\[
 x(t)= F^{-1} (\eta'(tq^{-1})).
\]
Podemos continuar usando esse mesmo argumento de reconstrução para uma extensão 
nos outros intervalos  $[aq^{n+1}, aq^n]$ com $n>2$, mas 
isso requer hipóteses ainda mais fortes: $\eta$ de classe $C^n$, 
compatibilidade das derivadas laterais nos extremos dos intervalos (para 
garantir continuidade) e que $F$ tenha inversa derivável.

%
%
%
%
%
%

\subsection{Equação linear e crescimento exponencial}

Muitas equações diferenciais lineares tem soluções formalmente como $\exp\{ 
\alpha t \} 
$, onde $\alpha$ é calculado a partir de uma equação característica. No caso de
uma equação linear com retardo 
constante $r>0$, $x'(t)= \lambda x(t-r)$, soluções deste tipo ainda existem, 
onde 
$\alpha$ são soluções de equações transcendentes do tipo $ \alpha = \lambda 
e^{-\alpha r}$. Isso garante a possibilidade de algum crescimento exponencial, 
um 
fenômeno típico de soluções de equações lineares. No entanto, as equações 
lineares com retardo proporcional $x'(t)= \lambda  
x(qt)$ não têm equações características (verifique!). Mostraremos que 
as soluções não têm comportamento exponencial, 
é o que mostra os próximos dois resultados:

\begin{teorema} A solução de $x'(t)= \lambda x(qt)$ com $q\in (0,1)$, $\lambda 
\geq 0$ e $x(0)= x_0 \in 
\R^*$ não tem crescimento exponencial, isto é
\[
 \lim_{t \rightarrow \infty} \ \frac{1}{t}\ \log |x(t)|=0.
\]
\end{teorema}

\bigskip

\noindent {\it Demonstração:} Usando o fato de que
\[
 \frac{d}{dt} \log |x(t)| = \frac{x'(t)}{x(t)} = \frac{\lambda x(qt)}{x(t)},
\]
temos que o crescimento exponencial é dado por
\[
 \lim_{t \rightarrow \infty} \ \frac{\lambda}{t}\ \int_{0}^t 
\frac{ x(qs)}{x(s)}\ ds.
\]
O resultado fica demonstrado se verificarmos que 
$
 \frac{x(qt)}{x(t)}
$
é decrescente e tende para zero. Mas de fato, para todo $t\geq 0$ (tomando 
$a=0$) temos
\[
 x(t)= x(qt) + \int_{qt}^t x(qs)\ ds.
\]
Assim, dividindo os dois lados da equação acima por $x(qt)$ e lembrando que 
$x(t)$ é crescente temos que o novo integrando  $\frac{x(qs)}{x(qt)} \geq 1$ 
para $s \in [qt, t]$. Portanto
$
\frac{x(t)}{x(qt)}
$ é crescente e tende a infinito quando $t$ vai para infinito.

\fdd

O limite no enunciado do teorema acima é conhecido na 
literatura como {\it expoente de Lyapunov}, neste caso esse expoente $\alpha 
=0$. É bem conhecido que em sistemas 
lineares canônicos (sem retardo) em dimensão superiores
esses expoentes dependem do subspaço onde está a condição inicial e são 
dados pela parte real dos autovalores da matriz de coeficientes. Voltando para 
o caso unidimensional e retardo pantográfico, mostramos que para coeficientes 
$\lambda<0$, tampouco temos comportamento 
exponencial:

\begin{teorema} A solução de $x'(t)= -\lambda x(qt)$ com $q\in (0,1)$, $\lambda 
> 0$ e $x(0)= x_0 \in 
\R^*$ muda de sinal.
\end{teorema}

\noindent {\it Demonstração:} Vamos mostrar inicialmente que a solução atinge  
o 
zero pelo menos uma vez no intervalo $\displaystyle \left(0, 
\frac{1}{\lambda(1-q)}\right]$. 
Sem perda 
de generalidade, vamos assumir que  $x(0)>0$. Temos que  $x(t)$ é sempre 
descrescente antes de chegar no zero. Suponha que $x(t)$ ainda seja  
estritamente positivo em  $[0, t_0]$, com $\displaystyle t_0 = 
\frac{q}{\lambda(1-q)}$ (caso contrário acabou essa parte da demonstração). 
Precisamos mostrar que neste caso  $x(t)$ cruza o zero no intervalo $(t_0, 
t_0q^{-1}]$. De fato, usando que para todo $s\in (t_0, t_0q^{-1}]$, vale  
$x(t_0)\leq x(qs)$, então temos, para todo $t$ neste mesmo intervalo:
\begin{eqnarray*}
 x(t) &=& x(t_0) -\lambda \int_{t_0}^{t} x(qs)\ ds \\
    && \\
& \leq & x(t_0) - \lambda (t-t_0) x(t_0) \\
      && \\
    & = & x(t_0) \big(1- \lambda (t-t_0)    \big) 
\end{eqnarray*}
Agora basta notar que o lado direito da desigualdade acima se anula quando 
$t=t_0 q^{-1}$.

Se chamarmos de $r_1$ a primeira 
raiz da solução, i.e. o tempo onde  $x(t)$  cruza 
o zero pela primeira vez, então no ponto de cruzamento temos que a solução 
ainda é decrescente, porque $x'(r_1)= -\lambda x(qr_1)<0$. Portanto a solução, 
depois de $r_1$ decresce até chegar em mínimo local estritamente negativo em  
$t=r_1 q^{-1}$, 
primeira vez onde a derivada se anula.

\fdd

\subsection{Comentários finais}

\noindent {\bf 1.} Ainda sobre a equação linear tratada no último teorema acima 
$x'(t)= -\lambda x(qt)$, com $\lambda >0$. Possivelmente exista uma demonstração 
elementar (como a demonstração do Teorema 2) de que a solução cruza o zero uma 
infinidade de vezes. Deixamos aqui essa questão para o leitor interessado. Outra 
pergunta que deixo em aberto neste contexto é a estabilidade desta solução: isto 
é, se a solução fica alternando valores positivos e negativos (não 
periodicamente),  os máximos e mínimos locais terão módulo limitado quando $t$ 
vai para infinito? Como essa estabilidade (ou não) depende dos parâmetros $q$ e 
$\lambda$? Para ver mais propriedades sobre o comportamento assintótico, neste caso e em outros mais gerais,
veja por exemplo \cite{Kato-McLeod}.

\bigskip

\noindent {\bf 2.} Outro fato curioso nas equações com retardo em geral é o 
aumento da regularidade da 
solução quando o tempo $t$ cresce. No caso de equação com retardo proporcional 
tratado 
aqui, note que tipicamente temos uma descontinuidade da derivada da solução no 
ponto $a$, da mesma maneira, uma descontinuidade da segunda derivada no ponto 
$aq^{-1}$ e assim por diante de modo que no ponto $aq^{-k}$ temos uma 
descontinuidade na $(k+1)$-ésima derivada. No interior dos intervalos 
$(aq^{-k}, aq^{-k-1})$ o mesmo fenômeno se passa: se a condição inicial tiver 
regularidade de classe $C^l$ no intervalo inicial
$[aq,a]$ então a solução vai ter regularidade
$C^{l+k}$ no intervalo $(aq^{-k+1}, a q^{-k})$. Ou em outras palavras, a 
solução $x(t)$ tem regularidade crescendo logaritmicamente: apresenta
\[
l + 1+ \lfloor  \frac{\log t - \log a}{\log q^{-1}} \rfloor
\]
derivadas contínuas em $t\in [qa, \infty) \setminus \{aq^{-k}; k=-1,0,1, \ldots 
\}$.

\bigskip

\noindent {\bf 3.} Sobre a propriedade de reconstrução (parágrafo 2.3 
acima). Pode ser interessante se tivermos explicitamente os critérios para a 
reconstrução 
da solução em intervalos ainda mais a esquerda, i.e. em $[q^ka, q^{k-1}a]$ para 
$k>2$. No seguinte sentido: em que condições poderíamos nos aproximar da origem 
$t=0$? E se isso for 
possível, em que medida essas condições não conduziriam para as soluções 
das 
equações mencionadas no parágrafo 2.2, onde a condição inicial degenera para 
um único ponto?
Ainda em outras palavras, será que essa reconstrução só é possível se a 
condição inicial $\eta$ for um fragmento de uma solução com condição inicial 
degenerada? Portanto, dado uma $F$, essa $\eta$ seria única (para cada valor 
limite em $t=0$)?

\bigskip

\noindent {\bf 4.} A exemplo das equações com retardo constante (ver por 
exemplo \cite{catuogno_ruffino} e as referências contidas ali), as equações 
pantográficas também fazem sentido em dinâmicas mais gerais. Por exemplo essas equações podem  aparecer em variedades diferenciáveis 
dotadas de uma conexão afim (riemannianas, por exemplo) que 
determine um transporte paralelo ao longo de 
curvas diferenciáveis.  De fato, dada 
uma curva  $\eta: [qa, T) \rightarrow M$ diferenciável, denote por $//_{s,t}: 
T_{\eta(s)}M \rightarrow T_{\eta(t)}M $  o transporte paralelo entre os 
pontos $\eta(s)$ e $\eta(t)$, para $qa 
\leq s \leq t$. Assim, uma 
equação com retardo proporcional nesta variedade  se escreve como
\[
x'(t)= //_{qt, t} X(x(qt)).
\]
com $X$ um campo de vetores, $a \leq t$ e com 
condições iniciais $x(t)= \eta (t)$ em $t\in [qa, 
a]$. Até onde nos consta, ainda estão em aberto questões relativas a 
propriedades das soluções das equações pantográficas neste 
contexto.



\end{document}